\documentclass[12pt, a4paper]{article}

\usepackage[usenames, dvipsnames]{color}
\usepackage[utf8]{inputenc}
\usepackage[english]{babel}
 
\usepackage{hyperref}
\hypersetup{
    colorlinks=true,
                   citecolor=red,
    linkcolor=blue, linktoc=all,
    filecolor=magenta,      
    urlcolor=blue, 
}
 
\usepackage{amsmath,amsxtra,amssymb,latexsym, amscd,amsthm}
\usepackage{indentfirst}
\usepackage{picinpar}
\usepackage{floatflt}
\usepackage{makeidx}   
\usepackage{longtable}%
\usepackage{multicol}%
\usepackage{enumerate}
\usepackage{hyperref}

\usepackage{tikz}
\usepackage{tikz-3dplot}
\usetikzlibrary{calc,fadings,decorations.pathreplacing}

\tikzset{%
  >=latex, 
  inner sep=0pt,%
  outer sep=2pt,%
  mark coordinate/.style={inner sep=0pt,outer sep=0pt,minimum size=3pt,
    fill=black,circle}%
}
\usepackage{tkz-euclide}
\usetikzlibrary{intersections}
\usetkzobj{all}

\usepackage{graphicx}

\theoremstyle{plain}

\def\g.#1.{\widehat{#1}}

\newcommand{\R}{\mathbb R}

\newtheorem{btdn}{Definition}
\newcommand{\dn}{\begin{btdn}\sl}
\newcommand{\hdn}{\end{btdn}}

\newtheorem{btvd}{Example}
\newcommand{\vd}{\begin{btvd}\sl}
\newcommand{\hvd}{\end{btvd}}

\newtheorem{btdl}{Theorem}[section]
\newcommand{\dl}{\begin{btdl}\sl}
\newcommand{\hdl}{\end{btdl}}

\newtheorem{btly}{Remark}[section]
\newcommand{\ly}{\begin{btly}}
\newcommand{\hly}{\end{btly}}

\newtheorem{btbd}{Lemma}[section]
\newcommand{\bd}{\begin{btbd}\sl}
\newcommand{\hbd}{\end{btbd}}

\newtheorem{btmd}{Proposition}[section]
\newcommand{\md}{\begin{btmd}\sl}
\newcommand{\hmd}{\end{btmd}}

\usepackage[left=2.3cm,right=2cm,top=3.8cm,bottom=3.1cm]{geometry}

\usepackage{fancyhdr}

\usepackage{wrapfig}
\usepackage{amsfonts}
{
\pagestyle{fancy}
\fancyhf{}

\pagestyle{fancy}
\fancyhf{}
{\chead{\sc {$q$-Moment\, Measures\, via\, Optimal\, Transport}}}
{\rhead{\thepage}}

\usepackage{tocloft}
\cftsetindents{section}{0em}{2em}
\cftsetindents{subsection}{0em}{2em}

\setcounter{tocdepth}{2}

\newcommand{\dd}{\mathrm{d}}

\newcommand{\res}{\mathop{\hbox{\vrule height 7pt width .5pt depth 0pt \vrule height .5pt width 6pt depth 0pt}}\nolimits}

\newcommand{\haus}{\mathcal H}

\newcommand{\Wc}{\mathcal{T}}

\newcommand{\ve}{\varepsilon}

\newcommand{\Pu}{\mathcal{P}_1(\R^d)}

\begin{document}

\textcolor[rgb]{1.0,1.0,1.0}{.}\\
\begin{center}
\noindent { {\large {\textbf{$\boldsymbol q$-MOMENT\, MEASURES\, AND\, APPLICATIONS:\\
\vspace{0.2cm}
A\, NEW\, APPROACH\, VIA\, OPTIMAL\, TRANSPORT}} }}\\
\vspace{0.7cm}
{{ { {\rm {\sc {Huynh}}\, \rm {\sc {Khanh}}\footnote{Institute of Mathematics, 
Vietnam Academy of Science and Technology, 18 Hoang Quoc Viet, Hanoi, Vietnam ({$\mathtt {khanh.edu02@gmail.com}$}).}}} \,\,\,\,\,\,{\rm {\sc {Filippo}}} \,{\rm {\sc {Santambrogio}}}\footnote{Institut Camille Jordan, Université Claude Bernard Lyon 1, 43 boulevard du 11 novembre 1918, 69622 Villeurbanne cedex, France, and Institut Universitaire de France ($\mathtt {santambrogio@math.univ}$-$\mathtt {lyon1.fr}$).}
} }

\renewenvironment {abstract}
{\begin{quote}
\noindent \rule{\linewidth}{.0pt}\par{\bfseries \textbf{\sc \abstractname}.}}
{\medskip
\end{quote}}
{\footnotesize {\begin{abstract}
In 2017, Bo’az Klartag obtained a new result in differential geometry on the existence of affine hemisphere of elliptic type. In his approach, a surface is associated with every a convex function $\varphi :{{\mathbb R}^n} \to \left( {0, + \infty } \right)$ and the condition for the surface to be an affine hemisphere involves the $2$-moment measure of $\varphi$ (a particular case of $q$-moment measures, i.e measures of the form ${\left( {\nabla \varphi } \right)_\# }{\varphi ^{ - \left( {n + q} \right)}}$ for $q > 0$). In Klartag's paper, $q$-moment measures are studied through a variational method requiring to minimize a functional among convex functions, which is studied using the Borell-Brascamp-Lieb inequality. In this paper, we attack the same problem through an optimal transport approach, since the convex function $\varphi$ is a Kantorovich potential (as already done for moment measures in a previous paper). The variational problem in this new approach becomes the minimization of a local functional and a transport cost among probability measures $\varrho$ and the optimizer $\varrho_{\rm {opt}}$ turns out to be of the form $\varrho_{\rm {opt}}  = {\varphi^{ - \left( {n + q} \right)}}$. \end{abstract}}}
\end{center}

{\footnotesize {\tableofcontents}}

\section{Introduction}
This paper stems for the analysis performed in \cite{K2,C1} and \cite{S3}. More precisely, Klartag studied in \cite {K2} the connection between the notion of affine hemisphere of elliptic type and that of $q$-{\it moment measures}, and proved an existence result thanks to a variational problem in the class of convex functions. To clarify what we mean, let us recall that, given a convex function $\varphi :{{\mathbb R}^n} \to \left( {0, + \infty } \right)$ with $\lim_{|x|\to\infty}\varphi(x)=+\infty$ and a positive real number $q > 0$, a Borel probability measure $\mu$ on ${\mathbb R}^n$ is said to be the $q$-moment measure of $\varphi$ if the differential of $\varphi$ pushes forward the measure ${\varphi ^{ - \left( {n + q} \right)}}{\rm d}x$ towards $\mu$:
\begin{align}
\label{e0}
\mu : = {\left( {\nabla \varphi } \right)_ {\ne} }\varrho ,\,\,\,\,\,\,\,\,{\rm {where}}\,\,{\rm d}\varrho  = {\varphi ^{ - \left( {n + q} \right)}}{\rm d}x.
\end{align}

A similar path was followed in \cite{C1} for the more well-known notion of {\it moment measures}, where \eqref{e0} is replaced by 
\begin{align}
\label{eMM}
\mu : = {\left( {\nabla \varphi } \right)_ {\ne} }\varrho ,\,\,\,\,\,\,\,\,{\rm {where}}\,\,{\rm d}\varrho  = {e^{-\varphi}}{\rm d}x.
\end{align}

Both \cite{C1} and \cite{K2} include a characterization of those measures $\mu$ which are indeed moment measures (or $q$-moment measures) of some function $\varphi$, by minimizing a convex functional on $\varphi$ involving $\mu$, and studying its subdifferential so as to prove that the minimizer solves \eqref{eMM} or \eqref{e0}. 

The experienced reader will have noticed the connection of these notions with that of optimal transport. Indeed, by Brenier's Theorem, the map $\nabla \varphi $ will be the optimal transport map for the quadratic cost $c\left( {x,y} \right) = \frac{1}{2}{\left| {x - y} \right|^2}$ from $\varrho$ to $\mu$. Based on this connection, the second author provided in \cite{S3} an optimal-transport-based approach to the problem of moment measures. The goal of the present paper is to do the same with $q$-moment measures.

Before doing this, let us clarify the connection of $q$-moment measures with affine spheres, which was the motivation of \cite{K2}. Indeed, Klartag obtained in \cite {K2} a new result on the existence and uniqueness of affine hemisphere of elliptic type (i.e, affine hemispheres whose center is inside a convex set bounded by a hypersurface) by studying the problem of $2$-moment measures. hHe showed that for any $q>1$ any probability measure $\mu$ with finite first moment, not supported on a hyperplane, and with $0$ as barycenter, can be represented as the $q$-moment measure for some (unique up to translations) essentially continuous convex function $\varphi :{{\mathbb R}^n} \to \left( {0, + \infty } \right)$. This can be applied to the problem of finding affine hemispheres of elliptic type if one chooses $\mu$ to be the uniform measure on a bounded convex set and takes $q=2$.

For the reader coming from optimal transport and not from convex geometry, it would be useful to explain the notion of affine spheres and hemispheres. Namely, given a smooth and connected hypersurface $M \subset {\mathbb R}^{n+1}$ which is locally strongly-convex and an arbitrary point $p \in M$, one considers the tangent space $P$ to $M$ at the point $p$ and its parallel translates ${\left( {{P_t}} \right)_{t \ge 0}}$ which are defined by ${P_0} = P$ and ${P_t} = {P_0} + tv$ where $v \notin P$  is a vector pointing to the convex side of $M$ at the point $p$. Then, for every $t$, the plane $P_t$ cuts the hypersurface and selects a compact sector $Q_t$ of finite volume inside $M$, and we call $G_t$ the barycenter of $Q_t$. As $t \to {0^ + }$ the point $G_t$ draws a curve which ends at $G_0 = p$. The tangent to the curve at $p$ selects a special direction, which is the direction of the affine normal. The line ${\ell_M}\left( p \right)$ passing through $p$ and directed as the affine normal at $p$ is called the affine normal line at $p \in M$. The hypersurface $M$ is called affinely-spherical of elliptic type with center at a point $x \in {\mathbb R}^{n+1}$ if all of the affine normal lines of $M$ meet at $x$ and the ray from an arbitrary point of $M$ to the center always passes through the convex side of it. Then, $M$ is called affine hemisphere of elliptic type with anchor $K$ if the following two conditions are satisfied: ${(i)}$ there exist compact, convex sets $K,\widetilde K   \subset {{\mathbb R}^{n + 1}}$ with ${\rm {dim}} \left( K \right) = n$ and ${\rm {dim}} \left( {\widetilde K  } \right) = n + 1$ such that $M$ does not intersect the affine hyperplane spanned by $K$ and $\partial \widetilde K   = M \cup K$;  ${(ii)}$ the center of $M$ belongs to the relative interior of $K$. All these notions are invariant for affine transformation (differently from the usual notion of normal to a surface), which explains the choice of the name.\\
\textcolor[rgb]{1.0,1.0,1.0}{......}
\begin{tikzpicture}[scale=0.8, join=round,decoration={markings,mark=between positions 0 and 1 step 1/4 with { \fill circle (1.5pt);}}, rotate=17.5]
\draw[line width=1](2.245,4.561)
.. controls (2.793,3.16) and (3.202,2.138) .. (3.612,1.473)
.. controls (4.5,0.018) and (5.434,0.154) .. (7.775,1.536)
.. controls (8.057,1.703) and (8.342,1.88) .. (8.668,2.082);
\draw [line width=1](2.246,4.561)
.. controls (2.492,5.015) and (4.147,4.649) .. (5.942,3.928)
.. controls (7.737,3.206) and (8.991,2.401) .. (8.744,2.131)
.. controls (8.496,1.861) and (6.84,2.227) .. (5.046,2.948)
.. controls (3.251,3.67) and (1.997,4.475) .. (2.315,4.610);
\draw[fill=blue!20!,line width=1](3.612,1.473)
.. controls (4.5,0.018) and (5.434,0.154) .. (7.775,1.536);
\draw[fill=blue!20!](3.612,1.473)
.. controls (3.449,1.753) and (3.545,2.027) .. (4.701,2.1)
.. controls (6.931,2.198) and (8.179,1.78) .. (7.775,1.536)
.. controls (7.371,1.293) and (6.112,1.051) .. (5.305,1.005)
.. controls (4.498,0.959) and (3.774,1.192) .. (3.612,1.473);
\draw (3.128,2.406)--(1.487,2.406)--(0,0)--(9.125,0)--(10.109,2.406)--(8.998,2.406);
\draw[dotted] (3.128,2.406)--(8.998,2.406);

\draw (7.8,-0.3)--(9.016,-0.3)--(10.000,2.106)--(8.978,2.106);   
\draw[dashed] (2,-0.3)--(7.8,-0.3);

\draw (7.8,-0.5)--(9.016,-0.5)--(10.000,1.906)--(8.678,1.906);   
\draw[dashed] (3,-0.5)--(7.8,-0.5);

\draw (7.8,-0.7)--(9.016,-0.7)--(10.000,1.706)--(8.678,1.706);   
\draw[dashed] (4,-0.7)--(7.8,-0.7);

\draw (7.8,-0.85)--(9.016,-0.85)--(10.000,1.556)--(8.678,1.556);   
\draw[dashed] (5,-0.85)--(7.8,-0.85);
\draw[postaction={decorate}] (5.071,0.503)
	.. controls (4.989,0.742) and (5.017,1.14) .. (5.145,1.378);

\draw[-stealth] (10.584,2.227)--(10.22,1.2); 

\draw[-stealth, line width=1.1pt] (5.071,0.503)--(4.248,3); 
\end{tikzpicture}\,\,\,\,\,\,
\tdplotsetmaincoords{60}{110}
\pgfmathsetmacro{\radius}{1}
\pgfmathsetmacro{\thetavec}{0}
\pgfmathsetmacro{\phivec}{0}
\begin{tikzpicture}[scale=3.1,tdplot_main_coords, rotate=-10]
\tdplotsetthetaplanecoords{\phivec}
\draw[dashed] (\radius,0,0) arc (0:360:\radius);
\shade[ball color=blue!10!white,opacity=0.2] (1cm,0) arc (0:-180:1cm and 5mm) arc (180:0:1cm and 1cm);
\draw (0, 1, 0) node [circle, fill=blue, inner sep=.02cm] () {};
\draw (0, 0, 1) node [circle, fill=green, inner sep=.02cm] () {};
\draw (-1, 0, 0) node [circle, fill=red, inner sep=.02cm] () {};
\tkzDefPoint(0,0){x}
\tkzDrawPoints[fill=black](x)
\tkzLabelPoints(x)
\draw[thick,->, line width=1.2] (0,0,1.15) -- (0,0,0.7) ;
\draw[dotted, line width=1.2] (0,0,0.7) -- (0,0,0) ;

\draw[thick,->, line width=1.2] (0,0.59,1.08) -- (0,0.3,0.549) ;
\draw[dotted, line width=1.2] (0,0.3,0.549) -- (0,0,0);

\draw[thick,->, line width=1.2] (0,-0.71,0.7156) -- (0,-0.2,0.201577) ;
\draw[dotted, line width=1.2] (0,-0.2,0.201577) -- (0,0,0);

\draw[thick,->, line width=1.2] (0,0.9,0.783) -- (0,0.25,0.2175) ;
\draw[dotted, line width=1.2] (0,0.25,0.2175) --(0,0,0);
\end{tikzpicture}\\
\textcolor[rgb]{1.0,1.0,1.0}{....}
{\scriptsize{ { {{\sc{Figure 1.}} Geometrical construction of the affine normal.}}}}\hspace{1.5cm}{\scriptsize{{\sc {Figure 2}}. An affine hemisphere of elliptic type in ${\mathbb R}^3$}}

 In \cite{K2}, Klartag proved that { {\sl{for any $n$-dimensional, compact and convex set ${K} \subset {\mathbb R}^{n+1}$, there always exists an $(n + 1)$-dimensional, compact and convex set ${\widetilde {K}} \subset {\mathbb R}^{n+1}$ whose boundary consists of two parts, the convex set $K$ itself is a facet and the rest of the boundary is an affine hemisphere with anchor ${K}$, which is centered at the Santaló point\footnote{Given any an $n$-dimensional, non-empty, bounded and convex set ${K} \subset {\mathbb R}^n$, then there exists a unique point $z$ in the interior of $K$ such that 
\begin{align}
\label{Santalo}
{\rm {Vol}}{_{\rm n}}\left( {{{\left( {{K} - z} \right)}^ \circ }} \right) = {\inf _{x \in {\mathop{\rm int}} \left( {K} \right)}}\left[ {{\rm {Vol}}{_{\rm n}}\left( {{{\left( {{K} - x} \right)}^ \circ }} \right)} \right]
\end{align}
here $S^{\circ}$ is the polar body of a convex set $S \subset {\mathbb R}^n$, defined by ${S^ {\circ} } = \left\{ {x \in {{\mathbb R}^n}:{{\sup }_{s \in S}}\left\langle {x,s} \right\rangle  \le 1} \right\}.$ The unique point $z$ in (\ref{Santalo}) is called the Santal\'o point of $K$ (see also \cite{A3, M2}). Note that the polar body ${{{\left( {{K} - z} \right)}^ \circ }}$ has its barycenter at the origin if and only if $z$ is the Santaló of $K$. As a corollary, the Santaló point of $K$ lies at the origin if and only if the barycenter of ${K}^{\circ}$ lies at the origin.} of $K$ and the affine hemisphere is uniquely determined up to transformations.}}} Thanks to affine transformations in ${\mathbb R}^{n+1}$, one may assume that the Santaló point of the set $K$ lies at the origin and 
$K \subset \left\{ {\left( {x,0} \right):x \in {{\mathbb R}^n}} \right\}.$ The set $K$ in those cases can be rewritten in form $K = {L^ \circ } \times \left\{ 0_{\mathbb R} \right\},$ where $L^{\circ}$ is the polar body of a certain convex set $L\subset  {{\mathbb R}^n}$. Then the proof of the existence of affine hemispheres of elliptic type relies upon analysis of the following PDE of Monge-Ampère type (see also \cite{F1})
\begin{align}
\label{e1}
\begin{cases}
\det {\nabla ^2}\varphi \left( x \right) \,\,\,&=\,\,\, \dfrac{{\rm {Vol}}_n\left(L\right)}{{{{\left( {\varphi \left( x \right)} \right)}^{n + 2}}}},\,\,\,\,\,\,\,\,x \in {{\mathbb R}^n}\\
\nabla \varphi \left( {{{\mathbb R}^n}} \right) &=\,\,\, L
\end{cases}
\end{align}
Indeed, it is proven (Theorem 1.2 in \cite{K2}) that, if $\varphi$ is a solution of (\ref{e1}) then we obtain an affine hemisphere with anchor $K = {L^ \circ } \times \left\{ 0_{\mathbb R} \right\}$ setting
\begin{align}
\label{e2}
M = \left\{ {\left( {\frac{x}{{\varphi \left( x \right)}},\,\,\frac{1}{{\varphi \left( x \right)}}} \right)} \in {{\mathbb R}^n} \times {\mathbb R}:\,\,x \in {\rm {dom}}\left( \varphi  \right)\right\}
\end{align}

Equation (\ref{e1}) is quite similar to the moment measures equation of Berman and Berndtsson \cite{1} in their work on Kähler-Einstein metrics in toric manifolds; Cordero-Erausquin and Klartag extended the study of \cite{1} in \cite{C1} presenting a functional version of the classical Minkowski problem or the logarithmic Minkowski problem and providing a variational characterization. More recently, the second author provided in \cite{S3} a dual counter-part of the results obtained in \cite{C1} with ideas coming from the theory of optimal transport, by considering the minimization of an entropy and a transport cost among probability measures. 

In \cite{K2}, the existence and the uniqueness of the solution to the equation (\ref{e1}) are proven via a variational method considering a minimization among convex functions. Namely, the author analyzes the subgradient of the functional
\begin{align}
\phi\mapsto {{\mathcal I}_{q}}\left( \phi \right) = {\left( {\int_{{{\mathbb R}^n}} {\frac{{{\rm d}x}}{{{{\left( {{\phi^*}\left( x \right)} \right)}^{n + q - 1}}}}} } \right)^{\frac{{ - 1}}{{q - 1}}}},
\end{align}
where ${\phi^*}\left( y \right): = {\sup _x}\left( {x.y - \phi\left( x \right)} \right)$ is the Legendre transform of $\phi$ and is the smallest function compatible with $\phi$ in the constraint $x.y \le \phi\left( x \right) + {\phi^*}\left( y \right).$ The convexity of the functional ${\mathcal I}_q$, which follows from the Borell-Brascamp-Lieb inequality, is the key for the proof. Equation \eqref{e1} is then obtained imposing that the uniform measure on $L$ belongs to the subdifferential, and taking $q=2$.

Taking ideas from the purely optimal-transport-based method for moment measures by the second author \cite{S3}, the main goal of the present work is to  reprove the same existence result of \cite{K2} with a different method, replacing functional inequalities techniques with ideas from optimal transport. In  \cite{S3} some heuristics are presented in order to guess which variational problem should be used (inspired by the theory of the JKO scheme for gradient flows, see for instance \cite{A1, JKO, S4}) and the functional to be minimized involves both a local functional\footnote{Local functionals over measures are defined as those functionals $F:{\mathcal P}\left( \Omega  \right) \to {\mathbb R}$ such that $F\left( {\mu  + \nu } \right) = F\left( \mu  \right) + F\left( \nu  \right)$ whenever  and $\mu$ and $\nu$ are mutually singular (i.e., there exists $A,B \subset \Omega $ with $A \cup B = \Omega $, $\mu \left( A \right) = 0$ and $\nu \left( B \right) = 0$).} of the form $\varrho\mapsto \int f(\varrho(x))\mathrm{d}x$ and a transport cost  ${{\mathcal T}}\left( {\varrho ,\mu } \right) $. We will not enter here into these considerations and we will directly look at the chosen functional.
We indeed study the following minimization problem 
\begin{align}
\left( P \right)\,\,\,\,\,\,\,\,\,\,\,\min \Big\{ {{\mathcal J}\left( \varrho  \right): = {\mathcal F}\left( \varrho  \right) + {\mathcal T}\left( {\varrho ,\mu } \right):\,\,\,\varrho  \in {{\mathcal P}_1}\left( {{{\mathbb R}^n}} \right)} \Big\},
\end{align}
where the functional ${\mathcal F}\left( . \right)$ is a very particular local functional (see Proposition 7.7 of \cite{S2}), given by 
\begin{equation}\label{defiF}
{\mathcal F}\left( \varrho  \right): = \int_{{{\mathbb R}^n}} {f\left( {{\varrho ^{{\rm {ac}}}}} \right){\rm d}x},\qquad \mbox{ writing }\varrho=\varrho^{{\rm ac}}(x){\rm d}x+\varrho^{{\rm sing}}\end{equation}
(the decomposition $\varrho=\varrho^{{\rm ac}}(x){\rm d}x+\varrho^{{\rm sing}}$ being the Radon-Nicodym decomposition of $\varrho$ into an absolutely continuous and a singular part w.r.t. the Lebesgue measure),
choosing  $f\left( t \right) =  - \frac{1}{\alpha }{t^\alpha }$ (to simplify notations, we set $\alpha  = 1 - \frac{1}{{n + q}} \in \left( {0,1} \right)$).
Note that the singular part of $\varrho$ does not appear in the functional since $f'(\infty)=0$.
The functional ${\mathcal T}\left( {.,\mu } \right)$ considered in this case is, instead, the maximal correlation functional (see \cite{S3}, section 3) and it is defined as follows 
\begin{align}
{\mathcal T}\left( {\varrho ,\mu } \right): = \sup \left\{ \int_{{\mathbb R}^n \times {\mathbb R}^n} {x.y} {\rm d}\gamma \left( {x,y} \right) \,\,\,|\,\,\,\,\gamma  \in \Gamma {\left( {\mu ,\varrho } \right)} \right\}
\end{align} 
where ${\Gamma {\left( {\mu ,\varrho } \right)}}$ 
is the set of transport plans between two probabily measures $\mu $ and $\varrho $ on ${\mathbb R}^n$. 

We will see in the rest of the paper that for some parts of the analysis we need $q>0$ (in particular for geodesic convexity issues) and for some other parts we need a stronger assumption, i.e. $q>1$ (in particular for the bounds which are needed to prove existence of a minimizer). However, the application to affine hemispheres (that we will not develop here, but is the motivation of the paper) requires $q=2$, which is fully covered by our results.

\section{Technical tools from optimal transport and local functionals on measures}
We recall here the main notions and notations that we will use throughout the paper. We refer to \cite{S2} (Chapters 1, 5 and 7) and to \cite{A1, V2, V3} for more details and complete proofs.\smallskip

\noindent\textbf{Optimal Transport and Wasserstein distances.} 

Given two probability measures $\mu ,\nu  \in {\mathcal P}\left( {{{\mathbb R}^n}} \right)$ we consider the set of transport plans
\begin{align}
\Gamma \left( {\mu ,\nu } \right) = \left\{ {\gamma  \in {\mathcal P}\left( {{{\mathbb R}^n} \times {{\mathbb R}^n}} \right):\,{{\left( {{\pi _x}} \right)}_\# }\gamma  = \mu ,\,{{\left( {{\pi _y}} \right)}_\# }\gamma  = \nu } \right\}
\end{align}
i.e. those probability measures on the product space having $\mu$ and $\nu$ as marginal measures.

For a cost function $c:{{\mathbb R}^n} \times {{\mathbb R}^n} \to \left[ {0, + \infty } \right]$ we consider the minimization problem
\begin{align}
\min \left\{ {\int {c{\rm d}\gamma } :\,\gamma  \in \Gamma \left( {\mu ,\nu } \right)} \right\}
\end{align}
which is called the Kantorovich optimal transport problem for the cost $c$ from $\mu$ to $\nu$. In particular, we consider the case $c\left( {x,y} \right) = \frac{1}{2}{\left| {x - y} \right|^2}$. In this case the above minimal value is finite whenever $\mu ,\nu  \in {{\mathcal P}_2}\left( {{{\mathbb R}^n}} \right)$, where ${{\mathcal P}_2}\left( {{{\mathbb R}^n}} \right): = \left\{ {\varrho  \in {\mathcal P}\left( {{{\mathbb R}^n}} \right):\,M_2(\varrho)<+\infty
} \right\}$ and for $p\geq 1$ we define $M_p(\varrho):=\int {{{\left| x \right|}^p}{\rm d}\varrho \left( x \right)}$.

For the above problem one can prove that the minimal value also equals the maximal value of a dual problem
\begin{align}
\max \left\{ {\int {\phi {\rm d}\mu }  + \int {\psi {\rm d}\nu } :\phi \left( x \right) + \psi \left( y \right) \le \frac{1}{2}{{\left| {x - y} \right|}^2}} \right\},
\end{align}
and that the optimal function $\phi$ may be used to construct an optimizer $\gamma$. Indeed, the optimal $\phi$ is locally Lipschitz and semiconcave (more precisely, $x \mapsto \frac{1}{2}{\left| x \right|^2} - \phi \left( x \right)$ is convex), and
differentiable $\mu$$-$a.e. if $\mu  \ll {{\mathcal L}^n}$; one can define a map $T:{{\mathbb R}^n} \to {{\mathbb R}^n}$ through $T\left( x \right) = x - \nabla \phi \left( x \right)$ and this map satisfies ${T_\# }\mu  = \nu $ ${\gamma _T}: = {\left( {{\rm {id}},T} \right)_\# }\mu $ (i.e. the image measure of $\mu$ through the map $x \mapsto \left( {x,T\left( x \right)} \right)$  belongs to ${\Gamma \left( {\mu ,\nu } \right)}$  and is optimal in the above problem). Moreover, the map $T$ is the gradient of the convex function $u$ given by $u\left( x \right) = \frac{1}{2}{\left| x \right|^2} - \phi \left( x \right)$  and is called the optimal transport map (for the quadratic cost $c\left( {x,y} \right) = \frac{1}{2}{\left| {x - y} \right|^2}$) from $\mu$ to $\nu$. The fact that the optimal transport map $T$ exists, is unique, and is the gradient of a convex function is known as Brenier Theorem (see \cite{2}).

The same could be obtained if one withdrew from the cost $\frac{1}{2}{\left| {x - y} \right|^2}$ the parts $\frac{1}{2}{\left| x \right|^2}$ and $\frac{1}{2}{\left| y \right|^2}$ which only depend on one variable each (hence, their integral w.r.t. $\gamma$ only depends on its marginals). Doing this we would get to the transport maximization problem
\begin{align}
{\mathcal T}\left( {\mu ,\nu } \right) = \max \left\{ {\int { {x.y}\, {\rm d}\gamma \left( {x,y} \right)\,\,\,\left|\,\,\, {\gamma  \in \Gamma \left( {\mu ,\nu } \right)} \right.} } \right\}
\end{align}
and the dual problem would become
\begin{align}
\inf \left\{ {\int {u{\rm d}\mu }  + \int {v{\rm d}\nu } :u\left( x \right) + v\left( x \right) \ge x.y} \right\}.
\end{align}
In this problem it is quite clear that any pair $(u, v)$ can be replaced with $\left( {u,{u^*}} \right)$ where ${u^*}\left( y \right) = {\sup _x}\left\{ {x.y - u\left( x \right)} \right\}$ is the Legendre transform of $u$, or even with $(u^{**},u^*)$. We can then assume that both $u$ and $v$ are convex and l.s.c.

 Then, it is easy to see by the primal-dual optimality conditions that the optimal $\gamma $ and the optimal $u$ satisfy
\begin{align}
{\rm {spt}}\left( \gamma  \right) \subset \left\{ {\left( {x,y} \right):u\left( x \right) + {u^*}\left( y \right) = x.y} \right\} = \left\{ {\left( {x,y} \right):y \in \partial u\left( x \right)} \right\},
\end{align}
which shows that $\gamma $ is concentrated on the graph of a map $T$ given by $T = \nabla u$, which is well-defined$\mu-$a.e provided $\mu  \ll {{\mathcal L}^n}$.

The value of the minimization problem with the quadratic cost may also be used to define a quantity, called Wasserstein distance, over ${{\mathcal P}_2}\left( {{{\mathbb R}^n}} \right)$
\begin{align}
{W_2}\left( {\mu ,\nu } \right): = {\left( {\min \left\{ {\int {{{\left| {x - y} \right|}^2}{\rm d}\gamma \left( {x,y} \right):\gamma  \in \Gamma \left( {\mu ,\nu } \right)} } \right\}} \right)^{1/2}}.
\end{align}
This quantity may be proven to be a distance over ${{\mathcal P}_2}\left( {{{\mathbb R}^n}} \right)$ and the space ${{\mathcal P}_2}\left( {{{\mathbb R}^n}} \right)$ endowed with the distance $W_2$ is called Wasserstein space of order 2, denoted by ${{\mathbb W}_2}\left( {{{\mathbb R}^n}} \right)$. On compact sets, this distance metrizes the usual weak convergence of probability measures (we say that a sequence $\varrho_n$ weakly converges to $\varrho$ if $\int \phi\,{\rm d}\varrho_n\to \int \phi\,{\rm d}\varrho$ for every bounded and continuous function $\phi$, and we write $\varrho_n\rightharpoonup\varrho$), while on $\R^n$ it metrizes a stronger notion of weak convergence, namely $\int \phi\,{\rm d}\varrho_n\to \int \phi\,{\rm d}\varrho$ for every continuous function $\phi$ with $|\phi(x)|\leq C(1+|x|^2)$.

The geodesics in this space play an important role in the theory of optimal transport. Given $\mu ,\nu  \in {{\mathcal P}_2}\left( {{{\mathbb R}^n}} \right)$ with $\mu  \ll {{\mathcal L}^n}$, we define ${\varrho _t}: = {\left( {\left( {1 - t} \right){\rm {id}} + tT} \right)_\# }\mu $, where $T$ is the optimal transport from $\mu$ to $\nu$. This curve ${\varrho}_t$ happens to be a constant speed geodesic for the distance $W_2$ connecting $\mu$ to $\nu$.

Once we know the geodesics in ${{\mathbb W}_2}\left( {{{\mathbb R}^n}} \right)$, one can wonder which functionals $F:{{\mathcal P}_2}\left( {{{\mathbb R}^n}} \right) \to {\mathbb R}$ are geodesically convex, i.e. convex along constant speed geodesics. This notion, applied to the
case of the Wasserstein spaces, is also called {\it displacement convexity} and has been introduced by McCann in \cite{M1}. It is very useful both to provide uniqueness results for variational problems and to provide sufficient optimality conditions. We will discuss this notion both in what co cerns the functional $\mathcal F$ and the functional $\mathcal T$.\newpage


\noindent\textbf{Maximal correlation functional.} 

We want to come back to the functional $\mathcal T$. It is useful to note that, if we have $\varrho, \mu  \in {{\mathcal P}_2}\left( {{{\mathbb R}^n}} \right)$, then we also have
\begin{align}
{\mathcal T}\left( {\varrho ,\mu } \right) = \frac{1}{2}\int_{{{\mathbb R}^n}} {{{\left| x \right|}^2}{\rm d}\varrho \left( x \right)}  + \frac{1}{2}\int_{{{\mathbb R}^n}} {{{\left| y \right|}^2}{\rm d}\mu \left( y \right)}  - \frac{1}{2}W_2^2\left( {\varrho ,\mu } \right).
\end{align}
If $\mathcal T$ comes from a transport cost, we may also observe that it stands for the maximal correlation between $\varrho$ and $\mu$, in the sense that we have
\begin{align}
{\mathcal T}\left( {\varrho ,\mu } \right) = \sup \Big\{ {{\mathbb E}\left[ {X.Y} \right]\,\,\,:\,\,\,X \sim \varrho ,Y \sim \mu } \Big\}.
\end{align}
For this reason, $\mathcal T$ will be called \textit{maximal correlation functional}. In this paper, we will make use of the following properties which were established in \cite{S3}.
\md
\label{p1}Suppose that $\mu\in\mathcal P(\R^n)$ is such that $M_1(\mu)<+\infty$ and $\int y,{\rm d}\mu(y)=0$. Then we have the following properties.\\
\,\,\,\,\,\,$\rm (1)$\,\, For every $\varrho  \in {{\mathcal P}_1}\left( {{{\mathbb R}^n}} \right),$ we have $0 \le {\mathcal T}\left( {\varrho ,\mu } \right) \le  + \infty .$\\
$\rm (2)$\,\, If $\varrho$ and ${\widetilde \varrho}$  are one obtained from one another by translation, then ${\mathcal T}\left( {\varrho ,\mu } \right) \,=\, {\mathcal T}\left( {\widetilde \varrho,\mu } \right).$\\
$\rm (3)$\,\, $\rm ($$ {Lower}$ $ {semicontinuity}$$\rm )$ If $\int {x{\rm d}{\varrho _n}\left( x \right)}  = 0$ and ${\varrho _n}$ $\rightharpoonup$ $\varrho$, then 
\begin{align}
{\mathcal T}\left( {\varrho ,\mu } \right) \,\,\,\le\,\,\, {\liminf _{n \to  + \infty }}\,{\mathcal T}\left( {{\varrho _n},\mu } \right).
\end{align}
$\rm (4)$\,\,There exists a sequence $\mu_n$ of compactly supported probability measures with ${\mu _n}$\,$\rightharpoonup$\,$\mu$ and $\int {y{\rm d}\mu \left( y \right)}  = 0$  such that for every $\varrho  \in {{\mathcal P}_1}\left( {{{\mathbb R}^n}} \right)$ we have ${\mathcal T}\left( {\varrho ,{\mu _n}} \right) \to {\mathcal T}\left( {\varrho ,\mu } \right)$.\\
$\rm (5)$\,\, If moreover $\mu$ is not supported on a hyperplane then, for every $\varrho$ such that $\int {x{\rm d}\varrho \left( x \right)}  = 0$, $\mathcal T$ satisfies an inequality of the form 
\begin{align}
{\mathcal T}\left( {\varrho ,\mu } \right) \,\,\,\ge\,\,\, c\int_{{{\mathbb R}^n}} {\left| x \right|{\rm d}\varrho \left( x \right)} 
\end{align}
 for $c = c\left( \mu  \right) > 0$, where
\[c\left( \mu  \right): = \frac{1}{{2n}}\inf \left\{ {\int_{{{\mathbb R}^n}} {\left| {y.e - \ell} \right|{\rm d}\mu \left( y \right):e \in {{\mathbb S}^{n - 1}},\ell \in {\mathbb R}} } \right\}.\]
$\rm (6)$\,\, $\rm ($$ {Displacement}$ $ {convexity}$$\rm )$. Let ${\varrho _0},{\varrho _1} \in {{\mathcal P}_2}\left( {{{\mathbb R}^n}} \right)$ be absolutely continuous measures, and let ${\varrho _t} = {\left( {\left( {1 - t} \right){\rm {id}} + tT} \right)_\# }{\varrho _0}$ be the unique constant speed geodesic connecting them for the Wasserstein distance $W_2$. Then $t \mapsto {\mathcal T}\left( {{\varrho _t},\mu } \right)$ is convex on $\left[ {0,1} \right]$. Moreover, if $\nabla v$ is the optimal transport from ${\varrho}_0$ to $\varrho_1$ and $\nabla \varphi$ is the optimal transport from ${\varrho}_0$ to $\mu$ then we have
\begin{align}
\dfrac{\rm d}{{{\rm d}t}}{\Big( {{\mathcal T}\left( {{\varrho _t},\mu } \right)} \Big)_{\left| {t = 0} \right.}} &\,\,\,\ge \displaystyle\int_{{{\mathbb R}^n}} {\left( {\nabla v\left( x \right) - x} \right).\nabla \varphi \left( x \right)\varrho_0 \left( x \right){\rm d}x}.
\end{align}
\hmd
\noindent\textbf{Lower semicontinuity and bounds for the local functional.} 
As we said, in this paper we will make use of a particular local functional $\mathcal F$, defined via \eqref{defiF}.

The functional $\mathcal F$ may be written in the form 
\begin{align}
\varrho  \mapsto \int_\Omega  {f\left( {{\varrho ^{{\rm {ac}}}}\left( x \right)} \right){\rm d}\lambda \left( x \right)}  \,+\, \left( {\mathop {\lim }\limits_{t \to  + \infty } \frac{{f\left( t \right)}}{t}} \right).{\varrho ^{{\rm {sing}}}}\left( \Omega  \right)
\end{align}
and for local functionals of this form lower semicontinuity results are well-known since at least \cite{B5} (see also Chapter 7 in \cite{S2} or Chapter 10 in \cite{A1}). Yet, these results are true  for convex and l.s.c. function $f:{{\mathbb R}_ + } \to {\mathbb R}$ (which is the case for $f\left( t \right) =  - \frac{1}{\alpha }{t^\alpha }$) but they require the reference measure $\lambda$ (which is the Lebesgue measure here) to be finite. Unfortunaltely, we need to consider here the whole, unbounded, space. Whenever $f\geq 0$ this difficulty can be easily solved by taking a sup over finite sub-measures of $\lambda$, but this is not the case here. We will then prove at the same time semicontinuity (in the spirit of Exercise 45 in \cite{S2}) and lower bounds on $\mathcal F$ in the following proposition. 
\md
Suppose $q>1$. Then, for each $\delta  \in \left( {\frac{n}{{n + q - 1}},1} \right)$, there exists a constant $C=C(\delta)$ such that the following estimate holds true
\begin{align}
{\mathcal F}\left( \varrho  \right) \,\,\,\ge\,\,\,  - C -M_1(\varrho)^\delta .
\end{align}
Moreover, whenever $\varrho_n$ is a sequence weakly-* converging to a measure $\varrho$ and $M_1(\varrho_n)$ is bounded, then we have $\mathcal F(\varrho)\leq \liminf_n\mathcal F(\varrho_n)$.
\hmd
\begin{proof}

First of all let us notice that when $f\left( t \right) =  - \frac{1}{\alpha }{t^\alpha }$, the Legendre transform $f^*$ of $f$, that is ${f^*}:{\mathbb R} \to {\mathbb R}$ defined by ${f^*}\left( h \right) = {\sup _{x \ge 0}}\left( {x.h - f\left( x \right)} \right) = {\sup _{x \ge 0}}\left( {x.h + \frac{1}{\alpha }{x^\alpha }} \right)$ is given by the following formula 
\begin{align}
{f^*}\left( h \right) \,=\,  \begin{cases}
 +\, \infty &{\rm {if}}\,\,\,h > 0,\\
\left( {\dfrac{1}{\alpha } - 1} \right){\left( { - h} \right)^{\dfrac{\alpha }{{\left( {\alpha  - 1} \right)}}}}&{\rm {if}}\,\,\,h \le 0.
\end{cases}
\end{align}

We first exploit the relation between $f$ and $f^*$ to obtain some lower bounds. We note that the functional $\mathcal F$ can be represented as follows
\begin{align}\label{ff*}
{\mathcal F}\left( \varrho  \right) &= \int_{{{\mathbb R}^n}} {\Big( {f\left( {{\varrho ^{{\rm {ac}}}}\left( x \right)} \right) + {f^*}\left( {h\left( x \right)} \right) - {\varrho ^{{\rm {ac}}}}\left( x \right)h\left( x \right)} \Big){\rm d}x} \notag\\
 &\,\,\,\,\,\,\,\,\,\,\,\,\,\,\,\,\,\,\,\,\,\,\,\,\,\,\,\,\,\,\,\,\,\,\,\,\,\,\,\,\,\,\,\,\,\,\,\,\,\,\,\,\,\,\,\,\,\,\,\,\,\,\,\,\,\,\,\,\,\,\,\,\,\,\,\,\,+ \int_{{{\mathbb R}^n}} {{\varrho ^{{\rm {ac}}}}\left( x \right)h\left( x \right){\rm d}x}  - \int_{{{\mathbb R}^n}} {{f^*}\left( {h\left( x \right)} \right){\rm d}x} 
\end{align}
for any function $h$ such that $h \in {L^1}\left( {{\varrho ^{{\rm {ac}}}}} \right)$ and ${f^*}\left( h \right) \in {L^1}\left( {{{\mathbb R}^n}} \right)$.

Assuming $M_1(\varrho)<+\infty$, we can fix $\delta<1$ and take $h\left( x \right) =  - {\left( {1 + \left| x \right|} \right)^\delta }$. This guarantees $h \in {L^1}\left( {{\varrho ^{{\rm {ac}}}}} \right)$; in order to have ${f^*}\left( h \right) \in {L^1}\left( {{{\mathbb R}^n}} \right)$ we need $\delta  \in \left( {\frac{n}{{n + q - 1}},1} \right)$. Indeed, using polar coordinates, we get
\begin{align}
0 \le \int_{{{\mathbb R}^n}} {{f^*}\left( {h\left( x \right)} \right){\rm d}x}  &= \left( {\frac{1}{\alpha } - 1} \right)\int_{{{\mathbb R}^n}} {{{\left( {1 + \left| x \right|} \right)}^{\frac{{ - \delta \alpha }}{{1 - \alpha }}}}{\rm d}x} \notag\\
 &= c\int_0^{ + \infty } {{r^{n - 1}}{{\left( {1 + r} \right)}^{\frac{{ - \delta \alpha }}{{1 - \alpha }}}}{\rm d}r} \notag
\end{align}
and the integral is finite as soon as 
\begin{align}
{n - 1 - }{\frac{{\delta \alpha }}{{1 - \alpha }}} <  - 1 \,\,\,\Longleftrightarrow \,\,\,\delta  > \frac{n}{{n + q - 1}}.
\end{align}

Let us start now from the first part of the claim, i.e. the lower bound. The Young-Fenchel inequality implies $f\left( {{\varrho ^{{\rm {ac}}}}\left( x \right)} \right) + {f^*}\left( {h\left( x \right)} \right) - {\varrho ^{{\rm {ac}}}}\left( x \right)h\left( x \right) \ge 0$, hence we get
\begin{align}
{\mathcal F}\left( \varrho  \right) \ge \int_{{{\mathbb R}^n}} {{\varrho ^{{\rm {ac}}}}\left( x \right)h\left( x \right){\rm d}x}  - \int_{{{\mathbb R}^n}} {{f^*}\left( {h\left( x \right)} \right){\rm d}x} .
\end{align}
The last term in this inequality is a constant independent of $\varrho$, and for the other term we use 
\begin{align}
\int_{{{\mathbb R}^n}} {{\varrho ^{{\rm {ac}}}}\left( x \right)h\left( x \right){\rm d}x}  &=  - \int_{{{\mathbb R}^n}} {{\varrho ^{{\rm {ac}}}}\left( x \right){{\left( {1 + \left| x \right|} \right)}^\delta }{\rm d}x} \notag\\
 &\ge   - \int_{{{\mathbb R}^n}} {{\varrho ^{{\rm {ac}}}}\left( x \right){\rm d}x}  - \int_{{{\mathbb R}^n}} {{\varrho ^{{\rm {ac}}}}\left( x \right){{\left| x \right|}^\delta }{\rm d}x} \notag\\
\end{align}
Applying H\"{o}lder inequality with exponents ${p_1} = \frac{1}{\delta }$ and ${p_2} = \frac{1}{{1 - \delta }}$, and using $ \varrho^{{\rm ac}}\leq \varrho$ and $\int \varrho^{{\rm ac}}\leq 1$, we have 
\begin{align}
\int_{{{\mathbb R}^n}} {{\varrho ^{{\rm {ac}}}}\left( x \right){{\left| x \right|}^\delta }{\rm d}x}  \,\,\,\le\,\,\, {\left( {\int_{{{\mathbb R}^n}} {\left| x \right|{\rm d}{\varrho }\left( x \right)} } \right)^\delta }.
\end{align}
This provides
\begin{align}
\int_{{{\mathbb R}^n}} {{\varrho ^{{\rm {ac}}}}\left( x \right)h\left( x \right){\rm d}x}  \ge  - 1 - {\left( {\int_{{{\mathbb R}^n}} {\left| x \right|{\rm d}\varrho \left( x \right)} } \right)^\delta }
\end{align}
Finally, for each $\delta  \in \left( {\frac{n}{{n + q - 1}},1} \right)$, the following estimate holds true
\begin{align}
{\mathcal F}\left( \varrho  \right) \,\,\,\ge\,\,\,  - C(\delta) - {\left( {\int_{{{\mathbb R}^n}} {\left| x \right|{\rm d}\varrho \left( x \right)} } \right)^\delta }.
\end{align}

We now have to prove the lower semicontinuity. We start from \eqref{ff*} which we re-write as 
\begin{align}\label{ff*}
{\mathcal F}\left( \varrho  \right) &= \int_{{{\mathbb R}^n}} {\Big( {f\left( {{\varrho ^{{\rm {ac}}}}\left( x \right)} \right) + {f^*}\left( {h\left( x \right)} \right) - {\varrho ^{{\rm {ac}}}}\left( x \right)h\left( x \right)} \Big){\rm d}x}-\int_{{{\mathbb R}^n}} h(x){\rm d}\varrho^{{\rm sing}}(x) \notag\\
 &\,\,\,\,\,\,\,\,\,\,\,\,\,\,\,\,\,\,\,\,\,\,\,\,\,\,\,\,\,\,\,\,\,\,\,\,\,\,\,\,\,\,\,\,\,\,\,\,\,\,\,\,\,\,\,\,\,\,\,\,\,\,\,\,\,\,\,\,\,\,\,\,\,\,\,\,\,+ \int_{{{\mathbb R}^n}} h(x) {\rm d}\varrho(x)  - \int_{{{\mathbb R}^n}} {{f^*}\left( {h\left( x \right)} \right){\rm d}x} .
\end{align}
Again, the function $h$ is fixed as $h(x)=(1+|x|)^\delta$, choosing $\delta$ as above.
We note that $\mathcal F$ is then composed of three parts. The last one is a constant, independent of $\varrho$. The previous one is 
$$\varrho\mapsto \int h {\rm d}\varrho.$$
We cannot say that this functional is continuous for the weak convergence of $\varrho$ since $h$ is not bounded. This is why we consider a sequence $\varrho_n$ with $M_1(\varrho_n)\leq C$, and we will exploit the sublinear behavior of $h$. If we fix an arbitrarily large constant $M$, we have
$$\left|\int h{\rm d}\varrho_n-\int h{\rm d}\varrho\right|\leq \left|\int \min\{h,M\}{\rm d}(\varrho_n-\varrho)\right|+\int_{\{x\,:\,h(x)>M\}}h\,{\rm d}(\varrho_n+\varrho).$$
Note that we have, for an arbitrary positive measure $\mu$ and $\delta<1$,
$$\int_{\{x\,:\,h(x)>M\}}h^{\frac 1 \delta}{\rm d}\mu \geq M^{\frac 1 \delta-1}\int_{\{x\,:\,h(x)>M\}}h\,{\rm d}\mu,$$
which, applied to $\mu=\varrho_n+\varrho$ and $h=(1+|x|)^\delta$, gives
$$2+2C\geq \int (1+|x|){\rm d}(\varrho_n+\varrho)\geq  M^{\frac 1 \delta-1}\int_{\{x\,:\,h(x)>M\}}h\,{\rm d}(\varrho_n+\varrho).$$
Hence we get 
$$\left|\int h{\rm d}\varrho_n-\int h{\rm d}\varrho\right|\leq \left|\int \min\{h,M\}{\rm d}(\varrho_n-\varrho)\right|+CM^{1-\frac 1 \delta}.$$
Taking the limsup in $n$ and using that $\min\{h,M\}$ is continuous and bounded we have 
$$\limsup_n\left|\int h{\rm d}\varrho_n-\int h{\rm d}\varrho\right|\leq CM^{1-\frac 1 \delta}.$$
Since $M$ is arbitrary, this provides 
$$\lim_n\int h{\rm d}\varrho_n=\int h{\rm d}\varrho.$$

Hence, this term in the expression of $\mathcal F$ is continuous for the convergence we use (weak convergence + bound on $M_1$). Finally, we just need to prove
$$\varrho\mapsto \int_{{{\mathbb R}^n}} {\Big( {f\left( {{\varrho ^{{\rm {ac}}}}\left( x \right)} \right) + {f^*}\left( {h\left( x \right)} \right) - {\varrho ^{{\rm {ac}}}}\left( x \right)h\left( x \right)} \Big){\rm d}x}-\int_{{{\mathbb R}^n}} h(x){\rm d}\varrho^{{\rm sing}}(x)$$
is l.s.c. for the weak convergence of probability measures. This functional is of the form
$$\int g(s,\varrho^{{\rm ac}}(x)){\rm d}x+ \int \lim_{t\to\infty}\frac{g(x,t)}{t}{\rm d}\varrho^{{\rm sing}}(x)$$
where $g(x,t)=f(t)+f^*(h(x))-h(x)t$. On a set of finite measure, this functional would be l.s.c. because of general results on local functionals (see, for instance, \cite{B5}). Here the reference measure has not finite mass, but the integrand $g$ is now positive, and hence the functional can be written as a supremum of local functionals of the same form on sets of finite measure, hence recovering lower semi-continuity.
\end{proof}

\noindent\textbf{Displacement convexity and strict displacement convexity}. 
Another important property of $\mathcal F$ concerns displacement convexity, the notion introduced by McCann in \cite{M1} and already presented in this paper concerning $\mathcal T$.  It is well-known from \cite{M1} that whenever $F(\varrho)$ has the the form $\int f(\varrho(x)){\rm d}x$, then displacement convexity is guaranteed as soon as $s\mapsto s^nf(s^{-n})$ is convex and decreasing, $n$ being the dimension of the ambient space. For the functional $\mathcal F$, this condition is satisfied as soon as $\alpha>1-\frac1n$, i.e. for $q>0$. We can now state the following proposition.

\md
\label{p2}
When restricted to ${\mathcal P}_2^{{\rm {ac}}}\left( {{{\mathbb R}^n}} \right)$, the functional $\mathcal F$ is displacement convex in ${\mathbb W}_2$ and strictly convex on every geodesic $t \mapsto {\varrho _t} = {\left( {\left( {1 - t} \right){\rm {id}} + tT} \right)_\# }\varrho $ unless the optimal map $T$ is such that $DT = I$ a.e. on $\{$${\varrho > 0}$$\}$. Moreover, if ${\varrho _t} = {\left( {\left( {1 - t} \right){\rm {id}} + tT} \right)_\# }\varrho $, then the derivative at $t = 0$ of $t \mapsto {\mathcal F}\left( {{\varrho _t}} \right)$  is given by 
\begin{align}
\label{di1}
\dfrac{\rm d}{{{\rm d}t}}{\Big( {{\mathcal F}\left( {{\varrho _t}} \right)} \Big)_{\left| {t = 0} \right.}} &\,\,\,= \left( {1 - \dfrac{1}{\alpha }} \right)\displaystyle\int_{{{\mathbb R}^n}} {{{\varrho  }^\alpha }{\rm {div}}\left( {T - {\rm {Id}}} \right){\rm d}x}.
\end{align}

In Formula (\ref{di1}), the divergence is to be taken in the a.e. sense, as $T$ is countably Lipschitz. Using $T = \nabla v$ with $v$ convex, this divergence is equal to ${\Delta ^{{\rm {ac}}}}v - n$, where ${\Delta ^{{\rm {ac}}}}$ is the absolutely continuous part of the distributional Laplacian of $v$, which is a positive measure.
\hmd

The proof of this statement follows the same line as that of points (3) and (4) in Proposition 2.1 of \cite{S3} (also see the computation in the Appendix A2 of \cite{BlaMosSan}).

\section{A Variational Principle for $\boldsymbol q$-Moment Measures}
As we sketched in the introduction, we consider the following variational problem. We fix $\mu  \in {{\mathcal P}_1}\left( {{{\mathbb R}^n}} \right)$ with $\int {y{\rm d}\mu \left( y \right)=0} $  and not supported on a hyperplane, and we want to solve
\begin{align*}
\left( P \right)\,\,\,\,\,\,\,\,\,\,\,\min \Big\{ {{\mathcal J}\left( \varrho  \right) = {\mathcal F}\left( \varrho  \right) + {\mathcal T}\left( {\varrho ,\mu } \right):\,\,\,\varrho  \in {{\mathcal P}_1}\left( {{{\mathbb R}^n}} \right)} \Big\}.
\end{align*}
\subsection{Existence of minimizers}
\dl
\label{the1}
Problem $(P)$ admits a solution.
\hdl
The idea for proving existence relies on the direct method in the calculus of variations, i.e. the fact that the minimized functional is lower semi-continuous and minimizing sequences are compact.

\textit{Proof of Theorem \ref{the1}.}\,\,\,Let ${\varrho _k} \in {{\mathcal P}_1}\left( {{{\mathbb R}^n}} \right)$ be a minimizing sequence, i.e. ${\inf _{\varrho  \in {{\mathcal P}_1}\left( {{{\mathbb R}^n}} \right)}}{\mathcal J}\left( \varrho  \right) = {\lim _{k \to  + \infty }}{\mathcal J}\left( {{\varrho _k}} \right).$ We can suppose that all $\varrho_k$ have $0$ as their barycenter as translations do not change the value of the two parts of the
functional and we may also assume that ${\mathcal J}\left( {{\varrho _k}} \right) \le {C_0}: = {\mathcal J}\left( {{\varrho ^0}} \right)$ with ${\varrho ^0} \in {{\mathcal P}_1}\left( {{{\mathbb R}^n}} \right)$. Using the estimates in Proposition \ref{p1}-part $(5)$ and Proposition \ref{p2}-part $(1)$, we have
\[{\mathcal J}\left( \varrho  \right) \ge {C_1}M_1(\varrho)-M_1(\varrho)^\delta - {C_2}.\]
This implies that $M_1(\varrho_k)$ must be bounded. According to Remark 5.1.5 in \cite{A1}, this gives tightness of the sequence $\varrho_k$ and hence we may extract a subsequence such that ${\varrho _{{k}}}$ $\rightharpoonup$ $\overline \varrho  $. On the other hand, by Proposition \ref{p2}-part $(2)$ we know that the functional $\mathcal F$ is l.s.c for the weak convergence when coupled to a bound on $M_1$, and the semi-continuity of $\mathcal T$ along sequences with $\int { x {\rm d}{\varrho _k}\left( x \right)}  = 0$ is in Proposition \ref{p1}-part $(3)$. These results imply 
\[{\mathcal J}\left( {\overline \varrho  } \right) \le \mathop {\lim \inf }\limits_{k \to  + \infty } {\mathcal J}\left( {{\varrho _k}} \right) = {\inf _{\varrho  \in {{\mathcal P}_1}\left( {{{\mathbb R}^n}} \right)}}{\mathcal J}\left( \varrho  \right)\]
which means that the minimum is indeed attained at ${\varrho _{{\rm {opt}}}} = \overline \varrho  .$  \hspace{5.35cm}$\square$
\ly 
Thanks to Theorem \ref{the4} in next section, each optimal solution to the problem $(P)$ is absolutely continuous. Hence, we can also obtain uniqueness of the solution from strict displacement convexity, using the geodesic $\varrho _t$ connecting two solutions $\varrho _0$ and $\varrho _1$ (see Proposition \ref{p2}-part (3); the case $DT = I$ a.e on $\left\{ {{\varrho _0} > 0} \right\}$ can be excluded exactly as in Lemma 4.2 in \cite{S3}).
\hly
\subsection{Properties of the optimal solutions}
The main goal of this section is to prove that each solution of the problem $(P)$ is absolutely continuous and show that if $\overline u$ is solution to the minimization problem
\begin{align}
\label{eqk}
{\rm {\inf}} \left\{ {\int_{{{\mathbb R}^n}} {u{\rm d}{{\varrho}_{\rm {opt}}} }  + \int_{{{\mathbb R}^n}} {{u^*}{\rm d}\mu } \,\,\,\left|\,\,\, {u:{{\mathbb R}^n} \to {\mathbb R} \cup \left\{ { + \infty } \right\}\,\,{\rm {convex\,\,and\,\,l.s.c}}} \right.} \right\}
\end{align}
 then there is a constant $c$ such that ${\overline u} > c$ and ${\varrho}_{\rm {opt}} = {{({\overline u} - c)}^{-{(n+q)}}}$. Finally, we will also prove that $\bar u$ is an essentially continuous convex function.
 
 To achieve our goals, we need some lemmas as follows.
\vspace{0.3cm}
\bd
\label{l1}
Given two densities $\sigma ,{\sigma _1} \in {L}^1\left( {{{\mathbb R}^n}} \right)$, with $\int_{{{\mathbb R}^n}} f(\sigma)>-\infty$, we have
\[\mathop {\lim }\limits_{\varepsilon  \to 0^+} \int_{{{\mathbb R}^n}} {\frac{{f\left( {\sigma  + \varepsilon \left( {{\sigma _1} - \sigma } \right)} \right) - f\left( \sigma  \right)}}{\varepsilon }{\rm d}x}  = \int_{{{\mathbb R}^n}} {f'\left( \sigma  \right)\left( {{\sigma _1} - \sigma } \right){\rm d}x} .\]
\hbd
\begin{proof}
By convexity of the function $f$, the inequality
\begin{align}
\label{eq3.1}
\int_{{{\mathbb R}^n}} {\frac{{f\left( {\sigma  + \varepsilon \left( {{\sigma _1} - \sigma } \right)} \right) - f\left( \sigma  \right)}}{\varepsilon }{\rm d}x}  \ge \int_{{{\mathbb R}^n}} {f'\left( \sigma  \right)\left( {{\sigma _1} - \sigma } \right){\rm d}x} 
\end{align}
is straightforward. Note that the right-hand side above could possibly be equal to $-\infty$. For the opposite inequality, we will use Fatou's Lemma. The pointwise convergence of the integrand is trivial and we can get an upper bound by means of the inequality
\[\frac{{f\left( {\sigma  + \varepsilon \left( {{\sigma _1} - \sigma } \right)} \right) - f\left( \sigma  \right)}}{\varepsilon } \le f\left( {{\sigma _1}} \right) - f\left( \sigma  \right)\leq - f\left( \sigma  \right),\]
valid for $\varepsilon  < 1$ (in the last inequality we used $f\le 0$). Then, since we supposed $f(\sigma)\in L^1$, we can apply Fatou's Lemma and we get
\begin{align}
\label{eq4.1}
\mathop {\lim \sup }\limits_{\varepsilon  \to 0} \int_{{{\mathbb R}^n}} {\frac{{f\left( {\sigma  + \varepsilon \left( {{\sigma _1} - \sigma } \right)} \right) - f\left( \sigma  \right)}}{\varepsilon }{\rm d}x}  \le \int_{{{\mathbb R}^n}} {f'\left( \sigma  \right)\left( {{\sigma _1} - \sigma } \right){\rm d}x} .
\end{align}
Combining (\ref{eq3.1}) and (\ref{eq4.1}), the lemma is proven.
\end{proof}

\bd
\label{l2.1}
If $\overline \varrho$ is a solution to the problem $(P),$ $\overline \varrho   = {\overline \varrho  ^{{\rm {ac}}}}{{\mathcal L}^n} + {\overline \varrho  ^{{\rm {sing}}}}$ and $\overline u :{{\mathbb R}^n} \to {\mathbb R} \cup \left\{ { + \infty } \right\}$ is solution to (\ref{eqk}), then
\begin{align}
\label{eq5.1}
\int_{{{\mathbb R}^n}} {\left( {\overline u  + f'\left( {{{\overline \varrho  }^{{\rm {ac}}}}} \right)} \right){\rm d}{{\overline \varrho  }^{{\rm {ac}}}}}  + \int_{{{\mathbb R}^n}} {\overline u {\rm d}{{\overline \varrho  }^{{\rm {sing}}}}}  \,\,\le\,\, \int_{{{\mathbb R}^n}} {\left( {\overline u  + f'\left( {{{\overline \varrho  }^{{\rm {ac}}}}} \right)} \right){\rm d}{\varrho ^{{\rm {ac}}}}}  + \int_{{{\mathbb R}^n}} {\overline u {\rm d}{\varrho ^{{\rm {sing}}}}} 
\end{align}
for all $\varrho  \in {{\mathcal P}_1}\left( {{{\mathbb R}^n}} \right)$ written as $\varrho   = {\varrho  ^{{\rm {ac}}}}{{\mathcal L}^n} + {\varrho  ^{{\rm {sing}}}}$.
\hbd
\begin{proof}
Assume that $\overline \varrho$ is optimal and $\overline u$ is a convex function realizing the minimum in the dual definition of ${\mathcal T}\left( {\overline \varrho  ,\mu } \right)$, then the functional
\[\varrho  \mapsto \int_{{{\mathbb R}^n}} {f\left( {{\varrho ^{{\rm {ac}}}}} \right){\rm d}x}  + \int_{{{\mathbb R}^n}} {\overline u {\rm d}\varrho }  = \int_{{{\mathbb R}^n}} {f\left( {{\varrho ^{{\rm {ac}}}}} \right){\rm d}x}  + \int_{{{\mathbb R}^n}} {\overline u {\rm d}{\varrho ^{{\rm {ac}}}}}  + \int_{{{\mathbb R}^n}} {\overline u {\rm d}{\varrho ^{{\rm {sing}}}}} \]
is minimal for ${\varrho} = {\overline \varrho}.$ \,Now for every $\varepsilon  \in \left( {0,1} \right)$ and any $\varrho  \in {{\mathcal P}_1}\left( {{{\mathbb R}^n}} \right)$,\, we define ${\varrho _\varepsilon }: = \left( {1 - \varepsilon } \right)\overline \varrho   + \varepsilon \varrho $, \,${\Phi _1}\left( {{\varrho ^{{\rm {ac}}}}} \right) = \int {f\left( {{\varrho ^{{\rm {ac}}}}} \right){\rm d}x} $,\, ${\Phi _2}\left( {{\varrho ^{{\rm {ac}}}}} \right) = \int {\overline u {\rm d}{\varrho ^{{\rm {ac}}}}} $ \,and\, ${\Phi _3}\left( {{\varrho ^{{\rm {sing}}}}} \right) = \int {\overline u {\rm d}{\varrho ^{{\rm {sing}}}}} .$ Using the optimality of $\overline \varrho$, we have
\begin{align}
\label{eq6.1}
\mathop {\lim }\limits_{\varepsilon  \to 0} \frac{{\Big\{ {\left( {{\Phi _1} + {\Phi _2}} \right)\left( {\varrho _\varepsilon ^{{\rm {ac}}}} \right) - \left( {{\Phi _1} + {\Phi _2}} \right)\left( {{\overline \varrho ^{{\rm {ac}}}}} \right)} \Big\} + \Big\{ {{\Phi _3}\left( {\varrho _\varepsilon ^{{\rm {sing}}}} \right) - {\Phi _3}\left( {{\overline \varrho ^{{\rm {sing}}}}} \right)} \Big\}}}{\varepsilon } \ge 0.
\end{align}
On the other hand,
\begin{align}
\label{eq7.1}
\begin{cases}
\dfrac{{{\Phi _2}\left( {\varrho _\varepsilon ^{{\rm {ac}}}} \right) - {\Phi _2}\left( {{{\overline \varrho  }^{{\rm {ac}}}}} \right)}}{\varepsilon } &= \displaystyle\int_{{{\mathbb R}^n}} {\overline u {\rm d}\left( {{\varrho ^{{\rm {ac}}}} - {{\overline \varrho  }^{{\rm {ac}}}}} \right)} \\
\dfrac{{{\Phi _3}\left( {\varrho _\varepsilon ^{{\rm {sing}}}} \right) - {\Phi _3}\left( {{{\overline \varrho  }^{{\rm {sing}}}}} \right)}}{\varepsilon } &= \displaystyle\int_{{{\mathbb R}^n}} {\overline u {\rm d}\left( {{\varrho ^{{\rm {sing}}}} - {{\overline \varrho  }^{{\rm {sing}}}}} \right)} 
\end{cases}
\end{align}
and by Lemma \ref{l1}, 
\begin{align}
\label{eq8.1}
\mathop {\lim }\limits_{\varepsilon  \to 0} \dfrac{{{\Phi _1}\left( {\varrho _\varepsilon ^{{\rm {ac}}}} \right) - {\Phi _1}\left( {{{\overline \varrho  }^{{\rm {ac}}}}} \right)}}{\varepsilon } = \displaystyle\int_{{{\mathbb R}^n}} {f'\left( {{{\overline \varrho  }^{{\rm {ac}}}}} \right){\rm d}\left( {{\varrho ^{{\rm {ac}}}} - {{\overline \varrho  }^{{\rm {ac}}}}} \right)} .
\end{align}
Combining (\ref{eq6.1}), (\ref{eq7.1}) and (\ref{eq8.1}), we obtain (\ref{eq5.1}). 
\end{proof}

\bd
\label{l3.1}
The solution $\overline \varrho$ to the problem $(P)$ cannot be such that ${{\overline \varrho}^{\rm {ac}}} = 0$ a.e.
\hbd
\begin{proof} Indeed, assume that ${{\overline \varrho}^{\rm {ac}}} = 0,$ then $\overline \varrho   = {\overline \varrho  ^{{\rm {sing}}}}$, ${\mathcal F}\left( {\overline \varrho  } \right) = {\mathcal F}\left( {{{\overline \varrho  }^{{\rm {sing}}}}} \right) = 0$ and hence 
\[{\mathcal J}\left( {\overline \varrho  } \right) = {\mathcal T}\left( {\overline \varrho  ,\mu } \right) \ge 0.\]
On the other hand, for $\eta  > 0,$ we take ${\varrho ^\eta } = \frac{{{{\boldsymbol 1}_{B\left( {0,\eta } \right)}}}}{{\left| {B\left( {0,\eta } \right)} \right|}}{{\mathcal L}^n}.$ Then ${\varrho ^\eta } \in {{\mathcal P}_1}\left( {{{\mathbb R}^n}} \right)$ and
\[{\mathcal T}\left( {{\varrho ^\eta },\mu } \right) = \max \left\{ {\int_{{\rm {spt}}\left( {{\varrho ^\eta }} \right) \times {{\mathbb R}^n}} {x.y{\rm d}\gamma :\gamma  \in \Gamma \left( {\mu ,{\varrho ^\eta }} \right)} } \right\} \le \eta \int_{{{\mathbb R}^n}} {\left| y \right|{\rm d}\mu }  = :C\eta \]
and value of the local functional ${\mathcal F}\left( {{\varrho ^\eta }} \right) =  - \frac{{\left| {B\left( {0,\eta } \right)} \right|}}{{{{\alpha\left| {B\left( {0,\eta } \right)} \right|}^\alpha }}} =  - {{{\alpha}^{-1}}\left| {B\left( {0,\eta } \right)} \right|^{1 - \alpha }} =  - {C_1}{\eta ^{n\left( {1 - \alpha } \right)}}.$ It follows that
\[{\mathcal J}\left( {{\varrho ^\eta }} \right) = C\eta  - {C_1}{\eta ^{n\left( {1 - \alpha } \right)}} .\]
Thanks to the assumption $\alpha>1-\frac 1n$ (corresponding to $q>0$), the above quantity can be made strictly negative as soon as $\eta$ is small enough, which gives a contradiction to the optimality of $\overline \varrho$.
\end{proof}

\dl
\label{the4}
If $\overline \varrho=\varrho_{\rm{opt}}$ is a solution to the problem $(P)$ and $\overline u$ is a solution to (\ref{eqk}), then there exists a constant $c$ such that $\overline u  \ge c$ and
\begin{align}
\label{eq8.2}
\left\{ \begin{array}{l}
\overline u  + f'\left( {{{\overline \varrho  }^{{\rm {ac}}}}} \right) = c\,\,\,{\rm {a.e}}\,\,{\rm {on}}\,\,\left\{ {{{\overline \varrho  }^{{\rm {ac}}}} > 0} \right\}\\
\\
\overline u  + f'\left( {{{\overline \varrho  }^{{\rm {ac}}}}} \right) \ge c\,\,\,{\rm {a.e}}\,\,{\rm {on}}\,\,\left\{ {{{\overline \varrho  }^{{\rm {ac}}}} = 0} \right\}
\end{array} \right.
\end{align}
and ${\overline \varrho}^{\rm {sing}}$ is concentrated on $\left\{ {\overline u  = c} \right\}$. Moreover, the set $\left\{ {\overline u  = c} \right\}$ is empty and thus ${\overline \varrho  ^{{\rm {sing}}}} = 0$, the optimal $\overline \varrho$ is absolutely continuous, with bounded density given by
\begin{align}
\overline \varrho   = {\overline \varrho  ^{{\rm {ac}}}} = \frac{1}{{{{\left( {\overline u  - c} \right)}^{n + q}}}}.
\end{align}
\hdl
\begin{proof}
Define the essential infimum 
\[c = {\rm {ess\,inf}}\Big\{ {\overline u  + f'\left( {{{\overline \varrho  }^{{\rm {ac}}}}} \right)} \Big\}: = \sup \Big\{ {\,\ell\,\,\,:\,\,\,{{\mathcal L}^n}\left( {\left\{ {x:\overline u \left( x \right) + f'\left( {{{\overline \varrho  }^{{\rm {ac}}}}\left( x \right)} \right) < \ell} \right\}} \right) = 0} \Big\}.\]
Then, we have $\overline u  + f'\left( {{{\overline \varrho  }^{{\rm {ac}}}}} \right) \ge c$ a.e, and so $\overline u  \ge \overline u  + f'\left( {{{\overline \varrho  }^{{\rm {ac}}}}} \right) \ge c$ a.e (since $f'\le 0$). Hence we have 
\begin{align}
\label{eq9.1}
\int_{{{\mathbb R}^n}} {\left( {\overline u  + f'\left( {{{\overline \varrho  }^{{\rm {ac}}}}} \right)} \right){\rm d}{{\overline \varrho  }^{{\rm {ac}}}}}  + \int_{{{\mathbb R}^n}} {\overline u {\rm d}{{\overline \varrho  }^{{\rm {sing}}}}}  \ge \int_{{{\mathbb R}^n}} {c{\rm d}{{\overline \varrho  }^{{\rm {ac}}}}}  + \int_{{{\mathbb R}^n}} {c{\rm d}{{\overline \varrho  }^{{\rm {sing}}}}}  = c.
\end{align}
Take any $c' > c$. By definition of essential infimum, the set $\left\{ {\overline u  + f'\left( {{{\overline \varrho  }^{{\rm {ac}}}}} \right) < c'} \right\}$ has positive Lebesgue measure, and so  we can choose $\varrho  \in {\mathcal P}_1^{{\rm {ac}}}\left( {{{\mathbb R}^n}} \right)$ such that $\varrho $ concentrated on the set $\left\{ {\overline u  + f'\left( {{{\overline \varrho  }^{{\rm {ac}}}}} \right) < c'} \right\}$. Then we get 
\begin{align}
\label{eq10.1}
\int_{{{\mathbb R}^n}} {\left( {\overline u  + f'\left( {{{\overline \varrho  }^{{\rm {ac}}}}} \right)} \right){\rm d}{\varrho ^{{\rm {ac}}}}}  + \int_{{{\mathbb R}^n}} {\overline u {\rm d}{\varrho ^{{\rm {sing}}}}}  = \int_{{{\mathbb R}^n}} {\left( {\overline u  + f'\left( {{{\overline \varrho  }^{{\rm {ac}}}}} \right)} \right){\rm d}\varrho }  < c'.
\end{align}
Combining (\ref{eq5.1}) in Lemma \ref{l2.1} and (\ref{eq10.1}), we obtain
\begin{align}
\label{eq11.1}
c' > \int_{{{\mathbb R}^n}} {\left( {\overline u  + f'\left( {{{\overline \varrho  }^{{\rm {ac}}}}} \right)} \right){\rm d}{{\overline \varrho  }^{{\rm {ac}}}}}  + \int_{{{\mathbb R}^n}} {\overline u {\rm d}{{\overline \varrho  }^{{\rm {sing}}}}} 
\end{align}
Letting $c' \to c$ in (\ref{eq11.1}), we get
\begin{align}
\label{eq12.1}
c \geq \int_{{{\mathbb R}^n}} {\left( {\overline u  + f'\left( {{{\overline \varrho  }^{{\rm {ac}}}}} \right)} \right){\rm d}{{\overline \varrho  }^{{\rm {ac}}}}}  + \int_{{{\mathbb R}^n}} {\overline u {\rm d}{{\overline \varrho  }^{{\rm {sing}}}}} 
\end{align}
From (\ref{eq9.1}) and (\ref{eq12.1}), we deduce that
\[\left\{ \begin{array}{l}
\displaystyle\int_{{{\mathbb R}^n}} {\left( {\overline u  + f'\left( {{{\overline \varrho  }^{{\rm {ac}}}}} \right)} \right){\rm d}{{\overline \varrho  }^{{\rm {ac}}}}}  = \int_{{{\mathbb R}^n}} {c{\rm d}{{\overline \varrho  }^{{\rm {ac}}}}} \\
\displaystyle\int_{{{\mathbb R}^n}} {\overline u {\rm d}{{\overline \varrho  }^{{\rm {sing}}}}}  \,\,\,\,\,\,\,\,\,\,\,\,\,\,\,\,\,\,\,\,\,\,\,\,\,\,\,= \int_{{{\mathbb R}^n}} {c{\rm d}{{\overline \varrho  }^{{\rm {sing}}}}} 
\end{array} \right.\]
It follows that 
\[\left\{ \begin{array}{l}
\overline u  + f'\left( {{{\overline \varrho  }^{{\rm {ac}}}}} \right) = c\,\,\,\,{\overline \varrho  ^{{\rm {ac}}}}{\rm {-a.e}}\\
\overline u  \,\,\,\,\,\,\,\,\,\,\,\,\,\,\,\,\,\,\,\,\,\,\,\,\,= c\,\,\,\,{\overline \varrho  ^{{\rm {sing}}}}{\rm {-a.e}}
\end{array} \right.\]
This implies that (\ref{eq8.2}) is satisfied and ${\overline \varrho}^{\rm {sing}}$ is concentrated on $\left\{ {\overline u  = c} \right\}$.

Setting $\varphi : = \overline u  - c$  we would like to prove that the set 
$\left\{ {\varphi  = 0} \right\}$ is empty. 

First, let us note that the interior of $\left\{ {\varphi  <  + \infty } \right\}$ is not empty. Indeed, should it be empty, 
since the set $\left\{ {\varphi  <  + \infty } \right\}$ is a convex set, then
it would be negligible, which means that in this case $\varphi  =  + \infty$ a.e. and hence ${\overline \varrho} = 0$, which is impossible. Thus the interior of $\left\{ {\varphi  <  + \infty } \right\}$ has to be non-empty. 
Next we suppose that there exists a point ${x_0} \in {\mathbb R}^n$ such that $\varphi\left(x_0\right) = 0$. 
First, we exclude the case where $x_0$ belongs to the interior of $\left\{ {\varphi  <  + \infty } \right\}$. Indeed, in this case there exists a neighborhood of $x_0$, denoted by 
${\mathcal N}\left( {{x_0}} \right) \Big(\subset \left\{ {\varphi  <  + \infty } \right\}\Big)$ 
where $\varphi $ is locally Lipschitz, so that we have
$\varphi  \leq L\left| {x - {x_0}} \right|$ for $x\in {\mathcal N}\left(x_0\right)$, which implies
\[1 > \int_{{\mathcal N}\left( {{x_0}} \right)} {\overline \varrho  {\rm d}x}  = \int_{{\mathcal N}\left( {{x_0}} \right)} {{\varphi ^{ - \left( {n + q} \right)}}{\rm d}x}  \ge \int_{{\mathcal N}\left( {{x_0}} \right)} {{{\Big( {L\left| {x - {x_0}} \right|} \Big)}^{ - \left( {n + q} \right)}}{\rm d}x}  =  + \infty .\]
This is impossible.

The case where $x_0$ lies on the boundary of the set $\left\{ {\varphi  <  + \infty } \right\}$ is more subtle. In this case we
choose $n$ points $x_1$, $x_2$, $...,$ $x_n$ in the interior of $\left\{ {\varphi  <  + \infty } \right\}$ 
so as to built a symplex $\Delta$ whose $\left(n + 1\right)$ vertices are $x_0$, $x_1$, $x_2$, $...$ and $x_n$, such that its interior is non-empty. 
On this symplex, the convex function $\varphi $ is finite and satisfies an inequality of the form
$\varphi  \le L\left| {x - {x_0}} \right|$ for all $x \in \Delta$, where the constant $L$ depends on $(x_i)_i$ and $(\varphi(x_i))_i$. Then, we find a similar contradiction as in the previous case.

We conclude that we cannot have $x_0 \in {\mathbb R}^n$ satisfying $\varphi\left(x_0\right) = 0.$ Since $\varphi$ is l.s.c. and we have $\lim_{|x|\to\infty}\varphi(x)=+\infty$ (this is a consequence of the integrability of $\bar\varrho=\varphi^{-(n+q)}$, then it admits a minimum on the space $\mathbb R^n$, and this minimum should be strictly positive. This proves the boundedness of $\bar\varrho$.
\end{proof}

We have now proven that the optimal $\bar \varrho$ can be expressed as $\bar \varrho=\varphi^{-(n+q)}$, and in order to fit the theory of Klartag we just need to prove that $\varphi$ is essentially continuous. This means that we want to prove $\lim_{x\to x_0} \varphi(x)=+\infty$ for $\haus^{d-1}$-a.e. $x_0\in \partial\{\varphi<+\infty\}$. 

In order to do this, given an optimal solution $\bar\varrho$, we choose a precise representative of it, and more precisely we take $\bar\varrho = \varphi^{-(n+q)}$, with $\varphi$ convex and l.s.c. We need to prove that $\bar\varrho$ vanishes on almost every point of the boundary.
\dl
Let $\bar\varrho$ be the precise representative above of a solution. Set $\Omega=\{\varphi<+\infty\}$. Then $\bar\varrho=0$ holds $\haus^{d-1}$-a.e. on $\partial\Omega$.
\hdl

\begin{proof} The proof will strongly follow that of Theorem 4.3 of \cite{S3}. 

Suppose $\bar\varrho>0$ on a set of positive $\haus^{d-1}$ measure on $\partial \Omega$. Writing locally $\Omega$ as $\{(x_1,x')\,:\, x_1>h(x')\}$ we assume that this set is given by $A=\{(x_1,x')\in\R\times\R^{d-1}\, x'\in B, x_1=h(x')\}$, where $B\subset \R^{d-1}$. Up to reducing the sets $A$ and $B$, we can suppose $\inf_{x\in A}\bar\varrho(x)>0$ for $x\in A$ and that $A,B$ are compact. We define $A_\ve:=\{(x_1,x')\in\R\times\R^{d-1}\, x'\in B, x_1\in[h(x'),h(x')+\ve]\}$. Then, by the continuity of $\bar\varrho$ inside $\Omega$, we can also assume, for small $\ve>0$, that we have $\bar\varrho(x_1,x')\geq c>0$ for $(x_1,x')\in A_\ve$. From the fact that $\bar\varrho$ is bounded, we also have the opposite inequality $\bar\varrho(A_\ve)\leq C_1\ve$.

Now, we define a new density $\varrho_\ve$ as a competitor by taking $T:A_\ve\to\R^d$ by $T(x_1,x')=(x_1-\ve,x')$ and setting
$$\varrho_\ve=\bar\varrho\res(A_\ve^c)+\frac 12 \bar\varrho\res(A_\ve)+\frac 12 T_\#(\bar\varrho\res(A_\ve)).$$
By computing the density of $\varrho_\ve$ we can check 
$$\mathcal F(\varrho_\ve)=\mathcal F(\bar\varrho)-(2^{\frac 1\alpha}-1)\int_{A_\ve}\bar\varrho\left(x\right)^\alpha{\rm d}x\leq \mathcal F(\bar\varrho)-C_2\ve.$$

In order to estimate $\Wc(\varrho_\ve,\mu)$, take the optimal function $\varphi$ (realizing $\Wc(\bar\varrho,\mu)=\int {\varphi}\,\dd\bar\varrho+\int {\varphi
}^*\,\dd\mu$) and modify it into a function $ \varphi_\delta$ as follows: take a convex, positive and superlinear function $\chi:\R^d\to\R$ with $\int \chi(x)\,\dd\mu(x)<+\infty$ (which exists because $\mu\in\Pu$), choose $\delta>0$ and set $ \varphi_\delta=( {\varphi}^*+\delta\chi)^*$. We have
$$ \varphi_\delta(x)\leq \varphi(x)\;\mbox{for every $x\in\Omega$}\quad  \varphi_\delta(T(x))\leq  \varphi(x)+\delta\chi^*\left(\frac{\ve e_1}{\delta}\right)\;\mbox{for every $x\in A_\ve$}$$
(where $e_1$ denotes te first vector of the canonical basis $e_1=(1,0,\dots,0)$).
Hence, 
$$\Wc(\varrho_\ve,\mu)\leq \int \varphi_\delta\,\dd\varrho_\ve+\int( \varphi^*+\delta\chi)\,\dd\mu\leq \Wc(\bar\varrho,\mu)+\frac 12\bar\varrho(A_\ve)\delta\chi^*\left(\frac{\ve e_1}{\delta}\right)+\delta\int\chi\,\dd\mu.$$
The optimality of $\bar\varrho$ compared to $\varrho_\ve$ provides
$$C_2\ve\leq \frac 12\bar\varrho(A_\ve)\delta\chi^*\left(\frac{\ve e_1}{\delta}\right)+\delta\int\chi\,\dd\mu.$$
Now, use $\bar\varrho(A_\ve)\leq C_1\ve$ and choose $\delta=c\ve$: we obtain, after dividing by $\ve$,
$$C_2\leq  \frac{C_1c}2\ve\delta\chi^*\left(\frac{e_1}{c}\right)+c\int\chi\,\dd\mu.$$
If we choose $c$ small enough, such that $c\int\chi\,\dd\mu<\frac 12 C_2$ we obtain a contradiction as $\ve\to 0$.
\end{proof}

\subsection{Sufficient optimality conditions}
To complete the current study, it remains to prove that every density $\varrho  = {\varphi ^{ - \left( {n + q} \right)}}$ such that ${\left( {\nabla \varphi } \right)_\# }\varrho  = \mu $ and $\varphi :{{\mathbb R}^n} \to \left( {0, + \infty } \right)$ is an essentially continuous convex function is necessarily a minimizer of ${\mathcal J}\left( \varrho  \right)$.  This would explain that the variational principle of the previous section finds exactly all the desired functions $\varphi$.

First of all, let us remind that for $\varrho$ to be integrable on $\R^n$ it is necessary that the convex function $\varphi$ satisfies $\lim_{|x|\to\infty} \varphi(x)=+\infty$ and that, because of convexity, the growth at infinity should be at least linear. Moreover, to be locally integrable around each point, the same computations as those that we showed at the end of the proof of Theorem \ref{the4} prove that $\varphi$ should be bounded from below by a strictly positive constant. In particular, $\varrho$ mus be a bounded density.

The result that we will prove is the following.
\dl
\label{dddd}
 Let $\varphi :{{\mathbb R}^n} \to \left( {0, + \infty } \right)$ be 
an essentially continuous convex function, consider $\overline \varrho   = {\varphi ^{ - \left( {n + q} \right)}}$ and suppose that 
$\mu  = \left( {\nabla \varphi } \right)_{\ne} \overline \varrho  .$ \,Then \,$\overline \varrho   \in {{\mathcal P}_1}\left( {{{\mathbb R}^n}} \right)$ 
\,and ${\mathcal J}\left( {\overline \varrho  } \right) = {\min _{\varrho  \in {{\mathcal P}_1}\left( {{{\mathbb R}^n}} \right)}}{\mathcal J}\left( \varrho  \right)$.
\hdl

Exactly as in the case of moment measures (\cite{S3}), the ideas to deal with the sufficient conditions come from displacement convexity (as it was the case in \cite{BlaMosSan}). Take an arbitrary $\varrho$ with compact support, and the geodesic curve ${\varrho _t} = {\left( {\left( {1 - t} \right){\rm {Id}} + tT} \right)_\# }\overline \varrho  $, 
where $T = \nabla v$ is the optimal transport from $\overline \varrho$ to $\varrho$. Then we know, from displacement convexity, that the following inequality holds true
\[{\mathcal J}\left( \varrho  \right) \,\,-\,\, {\mathcal J}\left( {\overline \varrho  } \right) \,\,\ge\,\, \frac{\rm d}{{{\rm d}t}}{\Big( {{\mathcal J}\left( {{\varrho _t}} \right)} \Big)_{\left| {t = 0} \right.}}\]
It is sufficient to show that the derivative in the right hand side above is non-negative. Using Proposition \ref{p1}-part $(6)$ and Proposition \ref{p2}-part $(4)$ we have
\begin{align}
\label{mmmm}
\begin{cases}
\dfrac{\rm d}{{{\rm d}t}}{\Big( {{\mathcal F}\left( {{\varrho _t}} \right)} \Big)_{\left| {t = 0} \right.}} &\,\,\,= \left( {1 - \dfrac{1}{\alpha }} \right)\displaystyle\int_{{{\mathbb R}^n}} {{\overline \varrho  }^\alpha {\rm {div}}\left( {T - {\rm {Id}}} \right){\rm d}x} \\
\dfrac{\rm d}{{{\rm d}t}}{\Big( {{\mathcal T}\left( {{\varrho _t},\mu } \right)} \Big)_{\left| {t = 0} \right.}} &\,\,\,\ge \displaystyle\int_{{{\mathbb R}^n}} {\left( {T\left( x \right) - x} \right).\nabla \varphi \left( x \right)\overline \varrho  \left( x \right){\rm d}x} .
\end{cases}
\end{align}
After considering $\varrho$ with compact support, we can use part $(5)$ of Proposition \ref{p2} to show that the optimality of $\bar\varrho$ is also valid when compared to non-compactly supported measures, by approximation.  

The arguments to be used to prove non-negativity of the derivative are very similar to the case of moment measures (\cite{S3}, Proposition 5.1), up to some modifications. The reader is invited to compare to \cite{S3}, and also to \cite{C1}. In particular, we first treat the terms involving the identity map in the above derivatives. This requires a computation similar to that in \cite{C1}, Lemma 5, i.e. the inequality
\begin{align}
\label{b1}
\frac{n}{{n + q - 1}}\int_{{{\mathbb R}^n}} {\frac{{{\rm d}x}}{{{\left(\varphi \left(x\right)\right) ^{n + q - 1}}}}}  \ge \int_{{{\mathbb R}^n}} {x.\nabla \varphi \left( x \right){\rm d}\overline \varrho  \left( x \right)} 
\end{align}
\noindent This is were we need the convex function $\varphi$ to be essentially continuous, as this guarantees $\displaystyle\int_{{{\mathbb R}^n}} {\nabla\left( {\frac{1}{{{\varphi ^{n + q - 1}}}}} \right)\left( x \right){\rm d}x}  = 0$ 
along the lines of (\cite{C1}, Lemma 4).
\bd
The inequality (\ref{b1}) holds for $\overline \varrho$ and $\varphi$ as in Theorem \ref{dddd}.
\hbd

\begin{proof} The computation is very similar to that developed for moment measures in \cite{C1}, Lemma 5. Thanks to the essentially-continuous property of the function $\varphi$, we have
\begin{align*}
\int_{{{\mathbb R}^n}} {\nabla \varphi \left( x \right){\rm d}\overline \varrho  \left( x \right)}  \,\,= \frac{{ - 1}}{{n + q - 1}}\int_{{{\mathbb R}^n}} {{\nabla}\left( {\frac{1}{{{\varphi ^{n + q - 1}}}}} \right)\left( x \right){\rm d}x}   \,\,=\,\, 0.
\end{align*}

Pick now a point $x_0$ in the interior of $\left\{\varphi < +\infty\right\}$ and define  $\mathcal K$ to be the class of all convex, smooth, compact sets $K$, contained in the interior of $\left\{\varphi < +\infty\right\}$ and containing  $x_0$ in their interior. Using the fact that the function $x\mapsto \nabla\varphi(x)\cdot (x-x_0)$ is bounded from below (by $\inf \varphi-\varphi(x_0)$), we get
\begin{align*}
&\int_{{{\mathbb R}^n}} {x.\nabla \varphi \left( x \right){\rm d}\overline \varrho  \left( x \right)}  \,\,\,=\,\,\, \int_{{{\mathbb R}^n}} {\nabla \varphi \left( x \right).\left( {x - {x_0}} \right){\rm d}\overline \varrho  \left( x \right)} 
 \leq \,\,\, \mathop {\sup }\limits_{K \in {\mathcal K}} \left\{ {\int_K {\nabla \varphi \left( x \right).\left( {x - {x_0}} \right){\rm d}\overline \varrho  \left( x \right)} } \right\}
 \end{align*}
 We then observe that, integrating by parts, we have
 $${\int_K {\nabla \varphi \left( x \right).\left( {x - {x_0}} \right){\rm d}\overline \varrho  \left( x \right)} }= {\frac{n}{{n + q - 1}}\int_K {\frac{{{\rm d}x}}{{{\varphi ^{n + q - 1}}}}}  - \frac{1}{{n + q - 1}}\displaystyle\int_{\partial K} {\frac{{\left( {x - {x_0}} \right).{{\boldsymbol {\rm n}}_x}}}{{{\varphi ^{n + q - 1}}}}{\rm d}\haus^{n-1}(x)} }.
$$
The second term of the right hand side is negative, using  the inequality $(x-x_0)\cdot{{\boldsymbol {\rm n}}_x}\geq 0$, valid for $K$ convex, $x_0\in K$ and $x\in \partial K$, and the first can be estimated by the integral on the whole space, which gives the claim. 
\end{proof}

\noindent\textit{Proof of Theorem \ref{dddd}.}\,\,\,Similar to Proposition 5.1 in \cite{S3}, using (\ref{mmmm}) and (\ref{b1}) above, we obtain
\begin{equation}\label{messadame}
{\mathcal J}\left( \varrho  \right) - {\mathcal J}\left( {\overline \varrho  } \right)\ge  {\frac{{ - 1}}{{n + q - 1}}\int_{\R^n} {\left( {{\Delta ^{{\rm {ac}}}}v} \right){{\overline \varrho  }}^\alpha }{\rm d}x}  + \int_{\R^n} {\nabla v\left( x \right).\nabla \varphi \left( x \right)\overline \varrho  \left( x \right){\rm d}x} .
\end{equation}
Note that we suppose that $\varrho$ is compactly supported, which implies $|\nabla v|\leq C$.  We will prove that the right hand side is positive, exactly as in  Proposition 5.1 in \cite{S3}, by first using $\Delta^{{\rm {ac}}} v\leq \Delta v$ (where $\Delta v $ is the distributional derivative of $v$). By integrating by parts we have 
$$\int_{B(0,R)} \Delta^{{\rm {ac}}} v\leq \int_{B(0,R)} \Delta v\leq \int_{\partial B(0,R)}|\nabla v|\leq CR^{n-1}.$$
Using $\varrho\leq cR^{-(n+q)}$ for large $R$, this implies $\int \Delta ^{{\rm {ac}}}v \overline \varrho <\infty$, since the integral on the annulus $B(0,2^k)\setminus B(0,2^{k-1})$ can be estimated by $C(2^k)^{n-1-(n+q)}$, which is summable. This shows that the first integral in \eqref{messadame} is well-defined and can be approximated with integrals on finite balls. To do the same for the second integral, we observe that we have $\nabla v\in L^\infty$ and $\nabla\varphi \in L^1(\bar\varrho)$ since this is equivalent to $M_1(\mu)<+\infty$, which gives $\nabla v\cdot\nabla\varphi \bar\varrho\in L^1(\mathbb R^n)$.

Hence, we have
\begin{align*}
&{\mathcal J}\left( \varrho  \right) - {\mathcal J}\left( {\overline \varrho  } \right)\\
 &\ge {\lim _{R \to  + \infty }}\Bigg( {\frac{{ - 1}}{{n + q - 1}}\int_{B\left( {0,R} \right)} {\left( {{\Delta }v} \right){{\overline \varrho  }^\alpha }{\rm d}x}  + \int_{B\left( {0,R} \right)} {\nabla v\left( x \right).\nabla \varphi \left( x \right)\overline \varrho  \left( x \right){\rm d}x} } \Bigg)\\
 &= {\lim _{R \to  + \infty }}\Bigg( {\frac{{ - 1}}{{n + q - 1}}\int_{B\left( {0,R} \right)} {\left( {{\Delta}v} \right){{\overline \varrho  }^\alpha }{\rm d}x}  + \int_{B\left( {0,R} \right)} {\nabla v\left( x \right).{\nabla _x}\left( {\frac{{ - 1}}{{n + q - 1}}.{\varphi ^{ - \left( {n + q - 1} \right)}}} \right)\left( x \right){\rm d}x} } \Bigg)\\
 &= \frac{{ - 1}}{{n + q - 1}}.{\lim _{R \to  + \infty }}\Bigg( {\int_{\partial B\left( {0,R} \right)} {\nabla v\left( x \right).{\boldsymbol {\rm n}}{{\Big( {\varphi \left( x \right)} \Big)}^{ - \left( {n + q - 1} \right)}}{\rm d}{{\mathcal H}^{n - 1}}\left( x \right)} } \Bigg) = 0,
\end{align*}
where the fact that the last integral tends to $0$ is justified using $|\nabla v|\leq C$ and bounding it from above by $CR^{n-1}R^{-(n+q-1)}$, which tends to $0$ thanks to $q>0$.
\begin{flushright}
$\square$
\end{flushright}

{\bf Acknowledgments} This work started as a part of the master thesis project of the first author, supervised by the second. His master studies at Universit\'e Paris-Saclay were partially supported by a scholarship from PGMO, a public grant part of the ``Investissement d'avenir'' project, reference ANR-11-LABX-0056-LMH, LabEx LMH, co-funded by the EDF company. The second author acknowledges the support of the Monge-Amp\`ere et G\'eom\'etrie Algorithmique project, funded by Agence nationale de la
recherche (ANR-16-CE40-0014 - MAGA).

{\footnotesize {

}}
\end{document}